\overfullrule=0pt
\centerline {\bf On the infimum of the upper envelope of certain families of functions}\par
\bigskip
\bigskip
\centerline {BIAGIO RICCERI}\par
\bigskip
\bigskip
{\it Abstract.} In this paper, given a topological space $X$, an interval $I\subseteq {\bf R}$ and five continuous functions $\varphi, \psi, \omega
:X\to {\bf R}$, $\alpha, \beta:I\to {\bf R}$, we are interested in the infimum of the function $\Phi:X\to ]-\infty,+\infty]$ defined by
$$\Phi(x)=\sup_{\lambda\in I}(\alpha(\lambda)\varphi(x)+\beta(\lambda)\psi(x))+\omega(x)\ .$$
Using a recent minimax theorem ([5]), we build a general scheme which provides the exact value of $\inf_X\Phi$ for a large class of functions $\Phi$. When additional compactness conditions are satisfied, our scheme provides also the existence of (explicitly detected) functions $\gamma, \eta:X\to {\bf R}$
such that, for some $\tilde x\in X$, one has
$$\gamma(\tilde x)\varphi(\tilde x)+\eta(\tilde x)\psi(\tilde x)+\omega(\tilde x)=\inf_{x\in X}(\gamma(\tilde x)\varphi(x)+\eta(\tilde x)\psi(x)+\omega(x))\ .$$
\bigskip
\bigskip
\bigskip
\bigskip
{\bf 1. Introduction}\par
\bigskip
Throughout the sequel, $X$ is a topological space, $I\subset {\bf R}$ is an interval and $\varphi, \psi, \omega:X\to {\bf R}$, $\alpha, \beta:I\to {\bf R}$ are five continuous functions.\par
\smallskip
 According to a standard terminology, the upper envelope of the family of functions $\{\alpha(\lambda)\varphi(\cdot)+
\beta(\lambda)\psi(\cdot)+\omega(\cdot)\}_{\lambda\in I}$ is the function $\Phi:X\to ]-\infty,+\infty]$ defined by
$$\Phi(x):=\sup_{\lambda\in I}(\alpha(\lambda)\varphi(x)+\beta(\lambda)\psi(x)+\omega(x))$$
for all $x\in X$.
The function $\Phi$ can have a quite complicated structure, depending on the specific choice of $\alpha$ and $\beta$. We are interested
in the infimum of $\Phi$. In this connection, we have the obvious estimate
$$\sup_{\lambda\in I}\inf_{x\in X}(\alpha(\lambda)\varphi(x)+\beta(\lambda)\psi(x)+\omega(x))\leq \inf_X\Phi\ .\eqno{(1.1)}$$
In general, the inequality in $(1.1)$ is strict. A particularly simple example
is provided by taking: $X=[0,2]$, $I=[0,1]$, $\varphi(x)=1-x^2$, $\psi(x)=x$, $\omega=0$, $\alpha(\lambda)=\lambda$, $\beta(\lambda)=1$.
Actually, for each $x\in [0,1]$, we have
$$\sup_{\lambda\in [0,1]}(\lambda(1-x^2)+x)=\cases{1-x^2+x & if $x\in [0,1]$\cr & \cr
x & if $x\in ]1,2]$\ .\cr}$$
Consequently
$$\inf_{[0,2]}\Phi=1\ .$$
Now, fix $\lambda\in [0,1]$. We have
$$\inf_{x\in [0,2]}(\lambda(1-x^2)+x)=\min\{\lambda, -3\lambda+2\}\ .$$
Therefore
$$\sup_{\lambda\in [0,1]}\inf_{x\in [0,2]}(\lambda(1-x^2)+x)={{1}\over {2}}\ ,$$
and so the inequality in $(1.1)$ is strict.\par
\smallskip
 Of course, it is interesting to know when equality occurs in $(1.1)$.\par
\smallskip
Before stating our contributions to this question, let us recall some definitions. Let $Y$ be a topological space. A function $\eta:Y\to {\bf R}$
is said to be inf-connected (resp. inf-compact) if, for each $r\in {\bf R}$, the set $\eta^{-1}(]-\infty,r[)$ 
(resp. $\eta^{-1}(]-\infty,r])$ is connected (resp. compact). If $-\eta$ is inf-compact, then $\eta$ is said to be
sup-compact.\par
\smallskip
The aim of this paper is to estabilish the following results whose proofs are based on the use of a quite recent minimax theorem, that is
Theorem 1.2 of [5].\par
\medskip
THEOREM 1.1. - {\it Let $I:=[a,b]$ be compact. Assume that:\par
\noindent
$(i_1)$\hskip 5pt there exists a set $D\subseteq I$, dense in $I$,  such that, for each $\lambda\in D$, the function
$\alpha(\lambda)\varphi(\cdot)+\beta(\lambda)\psi(\cdot)+\omega(\cdot)$ is inf-connected\ ;\par
\noindent
$(i_2)$\hskip 5pt the functions $\alpha, \beta$ are $C^1$ in $]a,b[$,
$\alpha'(\lambda)\neq 0$ for all $\lambda\in ]a,b[$ and the function $g:={{\beta'}\over {\alpha'}}$ is strictly monotone in $]a,b[$\ .\par
Finally, assume that one of the following conditions holds:\par
\noindent
$(i_3)$\hskip 5pt $g$ is increasing in $]a,b[$ and $\alpha'(\lambda)\psi(x)<0$ for all $(\lambda,x)\in ]a,b[\times (X\setminus \psi^{-1}(0))$\ ;\par
\noindent
$(i_4)$\hskip 5pt $g$ is decreasing in $]a,b[$ and $\alpha'(\lambda)\psi(x)>0$ for all $(\lambda,x)\in ]a,b[\times (X\setminus \psi^{-1}(0))$\ .\par
Set
$$\gamma:=\inf_{]a,b[}g\ ,$$
$$\delta:=\sup_{]a,b[}g$$
and let $h:{\bf R}\to [a,b]$ be the function defined by
$$h(\mu)=\cases{g^{-1}(\mu) & if $\mu\in g(]a,b[)$\cr & \cr
a & if $\mu\leq\gamma\hskip 3pt and\hskip 3pt (i_3)\hskip 3pt holds\hskip 3pt or\hskip 3pt \mu\geq\delta\hskip 3pt and\hskip 3pt
(i_4)\hskip 3pt holds$\cr & \cr
b & if  $\mu\leq\gamma\hskip 3pt and\hskip 3pt (i_4)\hskip 3pt holds\hskip 3pt or\hskip 3pt \mu\geq\delta\hskip 3pt and\hskip 3pt
(i_3)\hskip 3pt holds$\ .\cr}$$
Then,  we have
$$\Phi(x)=\cases{\alpha\left(h\left(-{{\varphi(x)}\over {\psi(x)}}\right)\right)\varphi(x)+
\beta\left(h\left(-{{\varphi(x)}\over {\psi(x)}}\right)\right)\psi(x)+\omega(x) & if $x\in X\setminus \psi^{-1}(0)$\cr & \cr
\max\{\alpha(a)\varphi(x),\alpha(b)\varphi(x)\}+\omega(x) & if $x\in\psi^{-1}(0)$\cr}$$
and
$$\inf_X\Phi=\sup_{\lambda\in I}\inf_{x\in X}(\alpha(\lambda)\varphi(x)+\beta(\lambda)\psi(x)+\omega(x))\ .$$}
\medskip
THEOREM 1.2. - {\it Assume that:\par
\noindent
$(i_0)$\hskip 5pt there exists a (possibly different) topology $\tau_1$ on $X$ such, for every $\lambda\in I$, the function
$\alpha(\lambda)\varphi(\cdot)+\beta(\lambda)\psi(\cdot)+\omega(\cdot)$ is $\tau_1$-lower semicontinuous and, for some $\lambda_0\in I$, the function $\alpha(\lambda_0)\varphi(\cdot)+\beta(\lambda_0)\psi(\cdot)+\omega(\cdot)$ is $\tau_1$-inf-compact\ ;\par
\noindent
$(i_1)$\hskip 5pt there exists a set $D\subseteq I$, dense in $I$,  such that, for each $\lambda\in D$, the function
$\alpha(\lambda)\varphi(\cdot)+\beta(\lambda)\psi(\cdot)+\omega(\cdot)$ is inf-connected\ ;\par
\noindent
$(i_2')$\hskip 5pt the functions $\alpha, \beta$ are derivable in $A$, where $\hbox {\rm int}(I)\subseteq A\subseteq I$,
$\alpha'(\lambda)\neq 0$ for all $\lambda\in A$, the function $g:={{\beta'}\over {\alpha'}}$ is strictly monotone in $A$ and $-{{\varphi(x)}\over {\psi(x)}}\in g(A)$ for all $x\in X\setminus \psi^{-1}(0)$.\par
Furthermore, assume that one of the following conditions holds:\par
\noindent
$(i_3')$\hskip 5pt $g$ is increasing in $A$ and $\alpha'(\lambda)\psi(x)<0$ for all $(\lambda,x)\in A\times (X\setminus \psi^{-1}(0))$\ ;\par
\noindent
$(i_4')$\hskip 5pt $g$ is decreasing in $A$ and $\alpha'(\lambda)\psi(x)>0$ for all $(\lambda,x)\in A\times (X\setminus \psi^{-1}(0))$\ .\par
Then, we have
$$\Phi(x)=\cases{\alpha\left(g^{-1}\left(-{{\varphi(x)}\over {\psi(x)}}\right)\right)\varphi(x)+
\beta\left(g^{-1}\left(-{{\varphi(x)}\over {\psi(x)}}\right)\right)\psi(x)+\omega(x) & if $x\in X\setminus \psi^{-1}(0)$\cr & \cr
\sup_{\lambda\in I}(\alpha(\lambda)\varphi(x)+\omega(x)) & if $x\in\psi^{-1}(0)$\cr}$$
and
$$\inf_X\Phi=\sup_{\lambda\in I}\inf_{x\in X}(\alpha(\lambda)\varphi(x)+\beta(\lambda)\psi(x)+\omega(x))\ .$$}
\medskip
THEOREM 1.3. - {\it Let the assumptions of Theorem 1.2 be satisfied. Moreover, assume that $\psi^{-1}(0)=\emptyset$ and that,
for some $x_0\in X$, the function
$\alpha(\cdot)\varphi(x_0)+\beta(\cdot)\psi(x_0)$ is sup-compact in $I$.\par
Then, there exists $\tilde x\in X$ such that
$$\alpha\left(g^{-1}\left(-{{\varphi(\tilde x)}\over {\psi(\tilde x)}}\right)\right)\varphi(\tilde x)+\beta\left(g^{-1}\left(-{{\varphi(\tilde x)}\over {\psi(\tilde x)}}\right)\right)\psi(\tilde x)+\omega(\tilde x)$$
$$=\inf_{x\in X}\left(\alpha\left(g^{-1}\left(-{{\varphi(\tilde x)}\over {\psi(\tilde x)}}\right)\right)\varphi(x)+
\beta\left(g^{-1}\left(-{{\varphi(\tilde x)}\over {\psi(\tilde x)}}\right)\right)\psi(x)+\omega(x)\right)\ . $$}
\medskip
An extremely important remark is that, even when, for each $\lambda\in I$, the function $\alpha(\lambda)\varphi(\cdot)+\beta(\lambda)\psi(\cdot)+\omega(\cdot)$ is convex, the above theorems cannot covered by the classical Sion's result ([6]), since, for some $x\in X$, the function
$\alpha(\cdot)\varphi(x)+\beta(\cdot)\psi(x)$ can be not quasi-concave.\par
\smallskip
The above results should properly be regarded as a general scheme from which one can derive, in a unitary way, boundless particular cases depending on the particular choices of the involved functions. Some samples will be highlighted in the third section. For instance, we will get the following
proposition whose statement is particularly simple:\par
\medskip
PROPOSITION 1.1. - {\it Let $X$ be a real Banach space, let $\varphi$ be linear and continuous and let $\psi$ be non-negative and Lipschitzian, with Lipschitz
constant equal to $\|\varphi\|_{X^*}$.\par
Then, for each constant $c\geq 1+\sqrt{2}$, one has
$$\inf_{x\in X}(c\varphi(x)+(c-\sqrt{2})\psi(x)+\sqrt{|\varphi(x)|^2+|\psi(x)|^2})=\left(c-{{1}\over {\sqrt{2}}}\right)\inf_{x\in X}(\varphi(x)+\psi(x))\ .$$}
\par
\bigskip
{\bf 2. Proofs of the basic results}\par
\bigskip
As we said, the main tool used to prove the above results is Theorem 1.2 of [5]. For the reader's convenience, here is its statement:\par
\medskip
THEOREM 2.A. - {\it Let $I$ be compact and let $f:X\times I\to {\bf R}$ be an upper semicontinuous
function such that $f(\cdot,\lambda)$ is continuous for all $\lambda\in I$. Assume also that, for each $\lambda$ in a dense subset of $I$, the function $f(\cdot,\lambda)$ is inf-connected and that,
for each $x\in X$, the set of all global maxima of the function $f(x,\cdot)$ is connected.\par
Then, one has
$$\sup_I\inf_Xf=\inf_X\sup_If\ .$$}
\bigskip
{\it Proof of Theorem 1.1}. First, notice that, since $g$ is continuous and strictly monotone, $g(]a,b[)$ agrees with the open interval
$]\gamma,\delta[$.  
Fix $x\in X\setminus \psi^{-1}(0)$. Let us show that the point $h\left(-{{\varphi(x)}\over {\psi(x)}}\right)$ is
the unique global maximum in $[a,b]$ of the function $\alpha(\cdot)\varphi(x)+\beta(\cdot)\psi(x)$. We prove this only for the case $(i_3)$. The
proof in the case $(i_4)$ is similar. So, assume that $(i_3)$ holds.
Assume that $-{{\varphi(x)}\over {\psi(x)}}\in ]\gamma,\delta[$.
Let $\tilde\lambda\in ]a,b[$ be such that
$${{\beta'(\tilde\lambda)}\over {\alpha'(\tilde\lambda)}}=-{{\varphi(x)}\over {\psi(x)}}\ .\eqno{(2.1)}$$
Since the function ${{\beta'}\over {\alpha'}}$ is increasing, we have
$${{\beta'(\lambda)}\over {\alpha'(\lambda)}}<-{{\varphi(x)}\over {\psi(x)}}$$
for all $\lambda\in ]a,\tilde\lambda[$ and
$${{\beta'(\lambda)}\over {\alpha'(\lambda)}}> -{{\varphi(x)}\over {\psi(x)}}$$
for all $\lambda\in ]\tilde\lambda,b[$. Then, by $(i_3)$, we have
$$\alpha'(\lambda)\varphi(x)+\beta'(\lambda)\psi(x)> 0$$
for all $\lambda\in ]a,\tilde\lambda[$ and
$$\alpha'(\lambda)\varphi(x)+\beta'(\lambda)\psi(x)< 0$$
for all $\lambda\in ]\tilde\lambda,b[$. Hence $\tilde\lambda$ is a global maximum of the function 
$\alpha(\cdot)\varphi(x)+\beta(\cdot)\psi(x)$. Finally,
since 
$$\alpha'(\lambda)\varphi(x)+\beta'(\lambda)\psi(x)\neq 0$$ 
for all $\lambda\in ]a,b[\setminus \{\tilde\lambda\}$, it follows that
$\tilde\lambda$ is the unique global maximum of that function in $[a,b]$ thanks to Rolle's theorem. From $(2.1)$ it follows that
$$\tilde\lambda=g^{-1}\left(-{{\varphi(x)}\over {\psi(x)}}\right)=
h\left(-{{\varphi(x)}\over {\psi(x)}}\right)\ .$$ Now, suppose that
$-{{\varphi(x)}\over {\psi(x)}}\not\in ]\gamma,\delta[$. If $-{{\varphi(x)}\over {\psi(x)}}\leq\delta$, then, by $(i_3)$, we have
$$\alpha'(\lambda)\varphi(x)+\beta'(\lambda)\psi(x)<0$$
for all $\lambda\in ]a,b[$. Hence, $a$ is the unique global maximum of the function $\alpha(\cdot)\varphi(x)+\beta(\cdot)\psi(x)$
in $[a,b]$. If $-{{\varphi(x)}\over {\psi(x)}}\geq\delta$, by $(i_3)$ again, we have
$$\alpha'(\lambda)\varphi(x)+\beta'(\lambda)\psi(x)>0$$
for all $\lambda\in ]a,b[$. So, $b$ is the unique global maximum of the function $\alpha(\cdot)\varphi(x)+\beta(\cdot)\psi(x)$
in $[a,b]$, and the claim is proved. Now suppose that $x\in \psi^{-1}(0)$. In this case, the set of all global maxima of the function
$\alpha(\cdot)\varphi(x)$ is either the whole $I$ (if $\varphi(x)=0$), or $\{a\}$, or $\{b\}$ in view of the strict monotonicity of
$\alpha$. Therefore, the function $f:X\times I\to {\bf R}$ defined by
$$f(x,\lambda)=\alpha(\lambda)\varphi(x)+\beta(\lambda)\psi(x)+\omega(x)$$
satisfies the hypotheses of Theorem 2.A, and so we have
$$\sup_{\lambda\in I}\inf_{x\in X}(\alpha(\lambda)\varphi(x)+\beta(\lambda)\psi(x)+\omega(x))=
\inf_{x\in X}\sup_{\lambda\in I}(\alpha(\lambda)\varphi(x)+\beta(\lambda)\psi(x)+\omega(x))\ .$$
In view of what seen above, the conclusion follows.\hfill $\bigtriangleup$\par
\bigskip
{\it Proof of Theorem 1.2.} Let $\{I_n\}$ be a non-decreasing sequence of compact intervals such that $I=\cup_{n\in {\bf N}}I_n$. Fix
$n\in {\bf N}$. Fix also $x\in X$. Denote by $M_n$ the set of all global maxima
of the restriction of the function $\alpha(\cdot)\varphi(x)+\beta(\cdot)\psi(x)$ to $I_n$. We claim that $M_n$ is connected. First, suppose
that $\psi(x)\neq 0$. If $-{{\varphi(x)}\over {\psi(x)}}\in g(A\cap I_n)$, then reasoning as in the proof of Theorem 1.1, we have
$M_n=\left\{g^{-1}\left(-{{\varphi(x)}\over {\psi(x)}}\right)\right\}$.
If $-{{\varphi(x)}\over {\psi(x)}}\not\in g(A\cap I_n)$, the function $\alpha(\cdot)\varphi(x)+\beta(\cdot)\psi(x)$ turns out to be strictly
monotone in $I_n$ (since its derivative does not vanish in int($I_n$)) and so we have either $M_n=\{\inf I_n\}$ or
$M_n=\{\sup I_n\}$. Now, let $\psi(x)=0$. In this case, we have $M_n=I_n$ if $\varphi(x)=0$. If $\varphi(x)\neq 0$, we still have either
$M_n=\{\inf I_n\}$ or $M_n=\{\sup I_n\}$, by strict monotonicity.
So, after noticing that $D\cap I_n$ is
dense in $I_n$, we realize that the function $f:X\times I_n\to {\bf R}$ defined by
$$f(x,\lambda)=\alpha(\lambda)\varphi(x)+\beta(\lambda)\psi(x)+\omega(x)$$
satisfies the hypotheses of Theorem 2.A. Consequently, we have
$$\sup_{\lambda\in I_n}\inf_{x\in X}(\alpha(\lambda)\varphi(x)+\beta(\lambda)\psi(x)+\omega(x))=
\inf_{x\in X}\sup_{\lambda\in I_n}(\alpha(\lambda)\varphi(x)+\beta(\lambda)\psi(x)+\omega(x))\ .$$
At this point, $(i_0)$ allows us to use Proposition 2.1 of [4] to obtain
$$\sup_{\lambda\in I}\inf_{x\in X}(\alpha(\lambda)\varphi(x)+\beta(\lambda)\psi(x)+\omega(x))=
\inf_{x\in X}\sup_{\lambda\in I}(\alpha(\lambda)\varphi(x)+\beta(\lambda)\psi(x)+\omega(x))\ .$$
Finally, in view of $(i_2')$, reasoning as in the proof of Theorem 1.1, we have
$$\sup_{\lambda\in I}(\alpha(\lambda)\varphi(x)+\beta(\lambda)\psi(x))=
\alpha\left(g^{-1}\left(-{{\varphi(x)}\over {\psi(x)}}\right)\right)\varphi(x)+
\beta\left(g^{-1}\left(-{{\varphi(x)}\over {\psi(x)}}\right)\right)\psi(x)$$
for all $x\in X\setminus \psi^{-1}(0)$, and we are done.\hfill $\bigtriangleup$\par
\bigskip
{\it Proof of Theorem 1.3.} Notice that the function $\sup_{\lambda\in I}(\alpha(\lambda)\varphi(\cdot)+\beta(\lambda)\psi(\cdot)+\omega(\cdot))$ 
turns out to be
$\tau_1$-lower semicontinous and $\tau_1$-inf-compact.
Consequently, there exists $\tilde x\in X$ such that
$$\inf_{x\in X}\sup_{\lambda\in I}(\alpha(\lambda)\varphi(\cdot)+\beta(\lambda)\psi(\cdot)+\omega(\cdot))= 
\sup_{\lambda\in I}(\alpha(\lambda)\varphi(\tilde x)+\beta(\lambda)\psi(\tilde x)+\omega(\tilde x))\ .\eqno{(2.2)}$$
Also, the function $\inf_{x\in X}(\alpha(\cdot)\varphi(x)+\beta(\cdot)\psi(x)+\omega(x))$ turns out to be sup-compact in $I$, and hence
there exists $\tilde\lambda\in I$ such that
$$\sup_{\lambda\in I}\inf_{x\in X}(\alpha(\lambda)\varphi(x)+\beta(\lambda)\psi(x)+\omega(x))=
\inf_{x\in X}(\alpha(\tilde\lambda)\varphi(x)+\beta(\tilde\lambda)\psi(x)+\omega(x))\ .\eqno{(2.3)}$$
Thanks to Theorem 1.2, the left-hand sides of $(2.2)$ and $(2.3)$ are equal. Therefore, we have
$$\alpha(\tilde\lambda)\varphi(\tilde x)+\beta(\tilde\lambda)\psi(\tilde x)+\omega(\tilde x)=\inf_{x\in X}(\alpha(\tilde\lambda)\varphi(x)+
\beta(\tilde\lambda)\psi(x)+\omega(x))=\sup_{\lambda\in I}(\alpha(\lambda)\varphi(\tilde x)+\beta(\lambda)\psi(\tilde x)+\omega(\tilde x))\ .$$
So, $\tilde\lambda$ is a global maximum of the function $\alpha(\cdot)\varphi(\tilde x)+\beta(\cdot)\psi(\tilde x)$ and hence, in view
of $(i_2')$, reasoning as in the proof of Theorem 1.1, we have 
$\tilde\lambda=g^{-1}\left(-{{\varphi(\tilde x)}\over {\psi(\tilde x)}}\right)$, and we are done.
\hfill $\bigtriangleup$
\bigskip
{\bf 3. Some applications}\par
\bigskip
The following proposition will be used in the subsequent alternative theorem.\par
\medskip
PROPOSITION 3.1. - {\it Let $X$ be a Hausdorff and let $f:X\to {\bf R}$ be a lower semicontinuous function such that, for
some $r>\inf_Xf$, the set $f^{-1}(]-\infty,r])$ is compact and disconnected.\par
Then, $f$ has at least two local minima.}\par
\smallskip
PROOF. Let $C_1, C_2$ be two non-empty, closed and disjoint subsets of $X$ such that 
$$f^{-1}(]-\infty,r])=C_1\cup C_2\ .$$
Notice that $C_1, C_2$ are compact. So, since $X$ is Hausdorff, there are two open and disjoint sets $\Omega_1, \Omega_2\subseteq X$ such that
$C_1\subseteq \Omega_1$, $C_2\subseteq \Omega_2$.
Since $f$ is lower semicontinuous, there are $x^*_1\in C_1$, $x^*_2\in C_2$ such that
$$f(x^*_1)=\inf_{C_1}f\ ,$$
$$f(x^*_2)=\inf_{C_2}f\ .$$
Without loss of generality, we assume that $f(x^*_1)\leq f(x^*_2)$. Fix $x\in \Omega_2$. We claim that $f(x^*_2)\leq f(x)$. Indeed, when
$f(x)>r$, the above inequality holds since $f(x^*_2)\leq r$. So, assume that $f(x)\leq r$. This implies that $x\in C_1\cup C_2$. Hence, since
$\Omega_2\cap C_1=\emptyset$, we have $x\in C_2$ and so the inequality is satisfied. Therefore, $x^*_2$ is a local minimum for $f$. But,
$x^*_1$ is even a global minimum for $f$, and the proof is complete.\hfill $\bigtriangleup$
\medskip
The following remarkable alternative theorem is a consequence of Theorem 1.3.\par
\medskip
THEOREM 3.1. - {\it Let $(i_2')$ be satisfied jointly with either $(i_3')$ or $(i_4')$. Moreover, assume that $X$ is Hausdorff and that,
for every $\lambda\in I$, the function $\alpha(\lambda)\varphi(\cdot)+\beta(\lambda)\psi(\cdot)+\omega(\cdot)$ is inf-compact. Finally,
suppose that $\psi^{-1}(0)=\emptyset$ and that, for some $x_0\in X$, the function $\alpha(\cdot)\varphi(x_0)+\beta(\cdot)\psi(x_0)$ is sup-compact in $I$.\par
Then, at least one of the following assertions holds:\par
\noindent
$(a)$\hskip 5pt there exists $\tilde x\in X$ such that
$$\alpha\left(h\left(-{{\varphi(\tilde x)}\over {\psi(\tilde x)}}\right)\right)\varphi(\tilde x)+\beta\left(h\left(-{{\varphi(\tilde x)}\over {\psi(\tilde x)}}\right)\right)\psi(\tilde x)+\omega(\tilde x)$$
$$=\inf_{x\in X}\left(\alpha\left(h\left(-{{\varphi(\tilde x)}\over {\psi(\tilde x)}}\right)\right)\varphi(x)+
\beta\left(h\left(-{{\varphi(\tilde x)}\over {\psi(\tilde x)}}\right)\right)\psi(x)+\omega(x)\right)\ ; $$
$(b)$\hskip 5pt there exists an open interval $J\subseteq I$ such that, for each $\lambda\in J$, the function
$\alpha(\lambda)\varphi(\cdot)+\beta(\lambda)\psi(\cdot)+\omega(\cdot)$ has at least two local minima.}\par
\smallskip
PROOF. Assume that $(a)$ does not hold. Notice that $(i_0)$ is satisfied if $\tau_1$ is the topology of $X$. Therefore,
 in view of Theorem 1.3, $(i_1)$ does not hold. That is, there exists an open interval
$J\subseteq I$ such that, for each $\lambda\in J$, some sublevel set of the function $\alpha(\lambda)\varphi(\cdot)+\beta(\lambda)\psi(\cdot)+\omega(\cdot)$ is not conncted. But, by assumption, such a sublevel set is also compact and hence the conclusion follows directly from Proposition
3.1.\hfill $\bigtriangleup$\par
\medskip
REMARK 3.1. - Theorem 3.1 can be usefully applied to get innovative multiplicity results for the solutions of classical nonlinear difference equations.\par
\medskip
Here is a consequence of Theorem 1.1.\par
\medskip
THEOREM 3.2. - {\it Let $\inf_X\psi\geq 0$ and let $c, d\in {\bf R}$ be such that the set
$$\left\{\lambda\in \left[-{{\pi}\over {2}},{{\pi}\over {2}}\right] : (\sin\lambda +c)\varphi(\cdot)+(\cos\lambda +d)\psi(\cdot)+\omega(\cdot)\hskip 3pt is\hskip 3pt
inf-connected\right\}$$
is dense in $\left[-{{\pi}\over {2}},{{\pi}\over {2}}\right]$.\par
Then, one has
$$\inf_{x\in X}(c\varphi(x)+d\psi(x)+\sqrt{|\varphi(x)|^2+|\psi(x)|^2}+\omega(x))
=\sup_{\lambda\in [-{{\pi}\over {2}},{{\pi}\over {2}}]}\inf_{x\in X}((\sin\lambda +c)\varphi(x)+(\cos\lambda +d)\psi(x)+\omega(x))\ .$$}
PROOF. We apply Theorem 1.1 taking 
$$I=\left[-{{\pi}\over {2}},{{\pi}\over {2}}\right]\ ,$$
 $$\alpha(\lambda)=\sin\lambda +c\ ,$$
$$\beta(\lambda)=\cos\lambda+d\ .$$
So that $$g(\lambda)=-\tan\lambda\ ,$$ $(i_4)$ holds and $$h(\mu)=-\arctan\mu$$ for all $\mu\in {\bf R}$.
Hence, if $x\in X\setminus \psi^{-1}(0)$, recalling that $\psi(x)>0$, we have
$$\Phi(x)=\left(\sin\arctan{{\varphi(x)}\over {\psi(x)}}+c\right)\varphi(x)+\left(\cos\arctan{{\varphi(x)}\over {\psi(x)}}+d\right)\psi(x)+\omega(x)$$
$$=\left({{\varphi(x)}\over {\sqrt{|\varphi(x)|^2+|\psi(x)|^2}}}+c\right)\varphi(x)+\left({{\psi(x)}\over {\sqrt{|\varphi(x)|^2+|\psi(x)|^2}}}+d\right)\psi(x)+\omega(x)$$
$$=c\varphi(x)+d\psi(x)+\sqrt{|\varphi(x)|^2+|\psi(x)|^2}+\omega(x)\ .$$
Now, the conclusion follows directly from Theorem 1.1 once we observe that
$$c\varphi(x)+d\psi(x)+\sqrt{|\varphi(x)|^2+|\psi(x)|^2}+\omega(x)=\max\{(c+1)\varphi(x), (c-1)\varphi(x)\}+\omega(x)$$
for all $x\in \psi^{-1}(0)$.\hfill $\bigtriangleup$
\medskip
In the same spirit as that of Theorem 3.2, we also obtain:\par
\medskip
THEOREM 3.3. - {\it Let $\psi(x)>0$ for all $x\in X$ and let $J\subseteq \left[-{{\pi}\over {2}},{{\pi}\over {2}}\right]$ be a compact interval such that
${{\varphi(x)}\over {\psi(x)}}\in \tan\left(J\cap \left]-{{\pi}\over {2}},{{\pi}\over {2}}\right[\right)$ for all $x\in X$. Let $c, d\in {\bf R}$.
Assume that the set
$$\left\{\lambda\in J : (\sin\lambda +c)\varphi(\cdot)+(\cos\lambda +d)\psi(\cdot)+\omega(\cdot)\hskip 3pt is\hskip 3pt
inf-connected\right\}$$
is dense in $J$ and that there exists
(a possibly different) topology $\tau_1$ on $X$ such that,
for every $\lambda\in J$, the function 
$(\sin\lambda +c)\varphi(\cdot)+(\cos\lambda +d)\psi(\cdot)+\omega(\cdot)$ is $\tau_1$-lower
semicontinuous and, for some $\lambda_0\in J$, the function 
$(\sin\lambda_0 +c)\varphi(\cdot)+(\cos\lambda_0 +d)\psi(\cdot)+\omega(\cdot)$ is $\tau_1$-inf-compact.
\par
Then, there exists $\tilde x\in X$ such that 
$$\left({{\varphi(\tilde x)}\over {\sqrt{|\varphi(\tilde x)|^2+|\psi(\tilde x)|^2}}}+c\right)\varphi(\tilde x)+
\left({{\psi(\tilde x)}\over {\sqrt{|\varphi(\tilde x)|^2+|\psi(\tilde x)|^2}}}+d\right)\psi(\tilde x)+\omega(\tilde x)$$
$$=\inf_{x\in X}\left(\left({{\varphi(\tilde x)}\over {\sqrt{|\varphi(\tilde x)|^2+|\psi(\tilde x)|^2}}}+c\right)\varphi(x)+
\left({{\psi(\tilde x)}\over {\sqrt{|\varphi(\tilde x)|^2+|\psi(\tilde x)|^2}}}+d\right)\psi(x)+\omega(x)\right)\ .$$}
PROOF. We apply Theorem 1.3 taking $I=J$, $A=J\cap \left]-{{\pi}\over {2}},{{\pi}\over {2}}\right[$ and proceed as in the proof of Theorem 3.2.\hfill $\bigtriangleup$\par
\medskip
Notice the following very particular corollary of Theorem 3.3:\par
\medskip
COROLLARY 3.1. - {\it Let $X$ be a closed convex subset of a reflexive real Banach space and let $\varphi$, $\psi$ be convex, with
$\varphi(x)\geq 0, \psi(x)>0$ for all $x\in X$. If $X$ is unbounded, suppose also that $\lim_{\|x\|\to +\infty}\psi(x)=+\infty$.\par
Then, there exists $\tilde x\in X$ such that
$$|\varphi(\tilde x)|^2+|\psi(\tilde x)|^2=\inf_{x\in X}(\varphi(\tilde x)\varphi(x)+\psi(\tilde x)\psi(x))\ .$$}
\smallskip
PROOF. Apply Theorem 3.3 taking $J=\left[0,{{\pi}\over {2}}\right]$ and $c=d=\omega=0$. The assumptions are satisfisfied since, 
 the function $(\sin\lambda)\varphi(\cdot)+(\cos\lambda)\psi(\cdot)$ is convex (and hence inf-connected) for all $\lambda\in J$ and
also weakly inf-compact if $\lambda<{{\pi}\over {2}}$. So, we can take the weak topology on $X$ as $\tau_1$. Then, by Theorem 3.3, there exists
$\tilde x\in X$ such that
$${{|\varphi(\tilde x)|^2}\over {\sqrt{|\varphi(\tilde x)|^2+|\psi(\tilde x)|^2}}}+
{{|\psi(\tilde x)|^2}\over {\sqrt{|\varphi(\tilde x)|^2+|\psi(\tilde x)|^2}}}
=\inf_{x\in X}\left({{\varphi(\tilde x)}\over {\sqrt{|\varphi(\tilde x)|^2+|\psi(\tilde x)|^2}}}\varphi(x)+
{{\psi(\tilde x)}\over {\sqrt{|\varphi(\tilde x)|^2+|\psi(\tilde x)|^2}}}\psi(x)\right)$$
and the conclusion follows multiplying both sides by $\sqrt{|\varphi(\tilde x)|^2+|\psi(\tilde x)|^2}$.\hfill $\bigtriangleup$\par
\medskip
In turn, from Corollary 3.1 we directly get\par
\medskip
COROLLARY 3.2. - {\it Let $X$ be a reflexive real Banach space and let $\varphi$, $\psi$ be convex, bounded below and Gateaux differentiable, with
 $\lim_{\|x\|\to +\infty}\psi(x)=+\infty$.\par
Then, for every $c\leq \inf_X\varphi$ and $d<\inf_X\psi$, there exists $\tilde x\in X$ such that
$$(\varphi(\tilde x)-c)\varphi'(\tilde x)+(\psi(\tilde x)-d)\psi'(\tilde x)=0\ .$$}\par
\medskip
A further corollary of Theorem 3.3, where $\varphi$ is not assumed to be convex, is as follows:\par
\medskip
COROLLARY 3.3. - {\it Let $X$ be a real Hilbert space and let $\varphi$ be $C^1$. Assume that $\varphi'$ is Lipschitzian, with Lipschitz constant
$L$. Moreover, let
$\eta:X\to ]0,+\infty[$ be a continuous, Gateaux differentiable, convex function such that
$$|\varphi(0)|+\|\varphi'(0)\|\|x\|\leq \eta(x) \eqno{(2.4)}$$
for all $x\in X$.\par
Then, there exists $\tilde x\in X$ such that
$$\varphi(\tilde x)\varphi'(\tilde x)+\left({{L}\over {2}}\|\tilde x\|^2+\eta(\tilde x)\right)(L\tilde x+\eta'(\tilde x))=0\ .$$}\par
\smallskip
PROOF. Since $\varphi$ is $C^1$, we have
$$\varphi(x)=\varphi(0)+\int_0^1\langle\varphi'(tx),x\rangle dt$$
for all $x\in X$. From this, by Lipschtzianity, we get
$$|\varphi(x)|\leq |\varphi(0)|+\int_0^1|\langle \varphi'(tx),x\rangle|dt\leq
|\varphi(0)|+\|x\|\int_0^1\|\varphi'(tx)\|dt\le\ |\varphi(0)|+\|x\|\int_0^1(\|\varphi'(0)\|+L\|x\|t)dt$$
$$=|\varphi(0)|+\|x\|\left(\|\varphi'(0)\|+{{L}\over {2}}\|x\|\right)$$
for all $x\in X$. Consequently, in view of $(2.4)$, we have
$${{|\varphi(x)|}\over {{{L}\over {2}}\|x\|^2+\eta(x)}}\leq 1$$
for all $x\in X$. Now, noticing that 
$$\sup_{\lambda\in \left [-{{\pi}\over {4}},{{\pi}\over {4}}\right ]}|\sin\lambda|=
\inf_{\lambda\in \left [-{{\pi}\over {4}},{{\pi}\over {4}}\right ]}\cos\lambda\ ,$$
for every $\lambda\in \left [-{{\pi}\over {4}},{{\pi}\over {4}}\right ]$, in view of a known result ([7], Corollary 42.7) the function
$(\sin\lambda)\varphi(\cdot)+{{L}\over {2}}(\cos\lambda)\|\cdot\|^2$ is convex. At this point,
we can apply Theorem 3.3 taking $$J=\left [-{{\pi}\over {4}},{{\pi}\over {4}}\right ]$$ and
$$\psi(x)={{L}\over {2}}\|x\|^2+\eta(x)\ .$$
The assumptions are satisfied with $c=d=\omega=0$. Indeed, for what above, the function $(\sin\lambda)\varphi(\cdot)+(\cos\lambda)\psi(\cdot)$ is
convex for all $\lambda\in \left [-{{\pi}\over {4}},{{\pi}\over {4}}\right ]$ and $${{\varphi(x)}\over {\psi(x)}}\in 
\tan\left(\left [-{{\pi}\over {4}},{{\pi}\over {4}}\right ]\right)=[-1,1]$$ for all $x\in X$. We then take the weak topology as $\tau_1$ and $\lambda_0=0$, and conclude as in the proof of Corollary 3.1.\hfill $\bigtriangleup$
\medskip
The next result is a further application of Theorem 1.1.\par
\medskip
THEOREM 3.4. - {\it  Let $I:=[a,b]$ be compact and assume that: $X$ is a real Banach space, $\varphi$ is continuous and linear, $\psi$ 
is Lipschitzian, with Lipschitz constant $L$. Moreover, assume that $(i_2)$ is satisfied jointly with either $(i_3)$ or $(i_4)$. Finally, assume
that the set
$$D:=\{\lambda\in I : |\beta(\lambda)|L<|\alpha(\lambda)|\|\varphi\|_{X^*}\}$$
is dense in $I$.\par
Then, if $h$ is defined as in Theorem 1.1, we have
$$\Phi(x)=\cases{\alpha\left(h\left(-{{\varphi(x)}\over {\psi(x)}}\right)\right)\varphi(x)+
\beta\left(h\left(-{{\varphi(x)}\over {\psi(x)}}\right)\right)\psi(x) & if $x\in X\setminus \psi^{-1}(0)$\cr & \cr
\max\{\alpha(a)\varphi(x),\alpha(b)\varphi(x)\} & if $x\in\psi^{-1}(0)$\cr}$$
and
$$\inf_X\Phi=\sup_{\lambda\in I\setminus D}\inf_{x\in X}(\alpha(\lambda)\varphi(x)+\beta(\lambda)\psi(x))\ .$$}
\smallskip
PROOF. Under our assumptions, by Theorem 2 of [3], for each $\lambda\in D$, the function $\alpha(\lambda)\varphi(\cdot)+\beta(\lambda)\psi(\cdot)$ is inf-connected and unbounded below. So
$$\sup_{\lambda\in I}\inf_{x\in X}(\alpha(\lambda)\varphi(x)+\beta(\lambda)\psi(x))=\sup_{\lambda\in I\setminus D}
\inf_{x\in X}(\alpha(\lambda)\varphi(x)+\beta(\lambda)\psi(x))$$
and the conclusion follows directly from Theorem 1.1.\hfill $\bigtriangleup$\par
\medskip
REMARK 3.2. - Other results in the same spirit as that of Theorem 3.4 can be found in [1], [2] and [3].\par
\medskip
Now, we are in a position to prove Proposition 1.1 stated in the Introduction.\par
\medskip
{\it Proof of Proposition 1.1.} Fix $c\geq 1+\sqrt{2}$. Observe that the minimum of the function $\lambda\to \sin\lambda-\cos\lambda$ in
$\left [-{{\pi}\over {2}},{{\pi}\over {2}}\right ]$ is equal to $-\sqrt{2}$ and is attained only at $-{{\pi}\over {4}}$. This clearly implies that, if
we take $$I=\left [-{{\pi}\over {2}},{{\pi}\over {2}}\right ]\ ,$$ $$\alpha(\lambda)=\sin\lambda +c$$ and $$\beta(\lambda)=\cos\lambda+c-\sqrt{2}\ ,$$
we have 
$$\{\lambda\in I : |\beta(\lambda)|<|\alpha(\lambda)|\}=I\setminus\left\{-{{\pi}\over {4}}\right\}\ . $$
Thus, since $(i_2)$ and $(i_4)$ hold, the assumptions of Theorem 3.4 are satisfied. So, proceedings as in the proof of Theorem 3.2, we get
$$\inf_{x\in X}(c\varphi(x)+(c-\sqrt{2})\psi(x)+\sqrt{|\varphi(x)|^2+|\psi(x)|^2})=\inf_{x\in X}\left(\alpha\left(-{{\pi}\over {4}}\right)\varphi(x)+
\beta\left(-{{\pi}\over {4}}\right)\psi(x)\right)$$
$$=\inf_{x\in X}\left(\left(c-{{1}\over {\sqrt{2}}}\right)\varphi(x)+\left(c-\sqrt{2}+{{1}\over {\sqrt{2}}}\right)\psi(x)\right)=
\left(c-{{1}\over {\sqrt{2}}}\right)\inf_{x\in X}(\varphi(x)+\psi(x))$$
and the proof is complete.\hfill $\bigtriangleup$\par
\medskip
The final result is a further consequence of Theorem 1.3.\par
\medskip
THEOREM 3.5. - {\it Let $\psi(x)>0$ for all $x\in X$, and assume that $\inf_X{{\varphi}\over {\psi}}>-\infty$. Set
$$I:=\left [\inf_X{{\varphi}\over {\psi}},\sup_X{{\varphi}\over {\psi}}\right ]\cap {\bf R}$$
and suppose that, for each $\lambda$ in a dense subset of $I$, the function $e^{\lambda}(\varphi(\cdot)+(1-\lambda)\psi(\cdot))+\omega(\cdot)$ is inf-connected.
Finally, assume that there exists (a possibly different) topology $\tau_1$ such that, for each $\lambda\in I$, the function
$e^{\lambda}(\varphi(\cdot)+(1-\lambda)\psi(\cdot))+\omega(\cdot)$ is $\tau_1$-lower semicontinuous
and, for some $\lambda_0\in I$, the function
$e^{\lambda_0}(\varphi(\cdot)+(1-\lambda_0)\psi(\cdot))+\omega(\cdot)$ is $\tau_1$-inf-compact.\par
Then, there exists $\tilde x\in X$ such that
$$\psi(\tilde x)+e^{-{{\varphi(\tilde x)}\over {\psi(\tilde x)}}}\omega(\tilde x)=
\inf_{x\in X}\left(\varphi(x)+\left(1-{{\varphi(\tilde x)}\over {\psi(\tilde x)}}\right)\psi(x)+e^{-{{\varphi(\tilde x)}\over {\psi(\tilde x)}}}\omega(x)\right)\ .$$}
\smallskip
PROOF. We are going to apply Theorem 1.3. In this connection, let $\alpha, \beta:I\to {\bf R}$ be defined by
$$\alpha(\lambda)=e^{\lambda}$$
and
$$\beta(\lambda)=(1-\lambda)e^{\lambda}$$
for all $\lambda\in I$. So that
$${{\beta'(\lambda)}\over {\alpha'(\lambda)}}=-\lambda\ .$$
Hence, $(i_2')$ and $(i_4')$ are satisfied, with $A=I$. Notice that, when $I$ is not bounded, we have
$$\lim_{\lambda\to +\infty}e^{\lambda}(\varphi(x)+(\lambda-1)\psi(x))=-\infty\ .$$
So, for each $x\in X$, the function $\lambda\to e^{\lambda}(\varphi(x)+(\lambda-1)\psi(x))+\omega(x)$ is sup-compact in $I$.
Therefore, all assumptions of Theorem 1.3 are satisfied. Consequently, there exists $\tilde x\in X$ such that
$$e^{{\varphi(\tilde x)}\over {\psi(\tilde x)}}\varphi(\tilde x)+\left(1-{{\varphi(\tilde x)}\over {\psi(\tilde x)}}\right)
e^{{\varphi(\tilde x)}\over {\gamma(\tilde x)}}\psi(\tilde x)+\omega(\tilde x)
=\inf_{x\in X}\left(e^{{\varphi(\tilde x)}\over {\psi(\tilde x)}}\varphi(x)+\left(1-{{\varphi(\tilde x)}\over {\psi(\tilde x)}}\right)
e^{{\varphi(\tilde x)}\over {\psi(\tilde x)}}\psi(x)+\omega(x)\right)$$
and hence, multiplying by $e^{-{{\varphi(\tilde x)}\over {\psi(\tilde x)}}}$, we get
$$\psi(\tilde x)+e^{-{{\varphi(\tilde x)}\over {\psi(\tilde x)}}}\omega(\tilde x)=
\inf_{x\in X}\left(\varphi(x)+\left(1-{{\varphi(\tilde x)}\over {\psi(\tilde x)}}\right)\psi(x)+e^{-{{\varphi(\tilde x)}\over {\psi(\tilde x)}}}\omega(x)\right)\ ,$$
as claimed.\hfill $\bigtriangleup$\par
\medskip
COROLLARY 3.4. - {\it Let the assumptions of Theorem 3.5 be satisfied. In addition, assume that $X$ is an open set in a real normed space and that
$\varphi, \psi$ are Gateaux differentiable.\par
Then, there exists $\tilde x\in X$ such that
$$\varphi'(\tilde x)+\left(1-{{\varphi(\tilde x)}\over {\psi(\tilde x)}}\right)\psi'(\tilde x)
+e^{-{{\varphi(\tilde x)}\over {\psi(\tilde x)}}}\omega'(\tilde x)=0\ .$$}
\medskip
REMARK 3.3. - Corollaries 3.2, 3.3 and 3.4 can be usefully applied to get existence theorems for different classes of nonlocal differential equations.
\bigskip
\bigskip
\bigskip
{\bf Acknowledgements:} This work has been funded by the European Union - NextGenerationEU Mission 4 - Component 2 - Investment 1.1 under the Italian Ministry of University and Research (MUR) programme "PRIN 2022" - grant number 2022BCFHN2 - Advanced theoretical aspects in PDEs and their applications - CUP: E53D23005650006. The author has also been supported by the Gruppo Nazionale per l'Analisi Matematica, la Probabilit\`a e 
le loro Applicazioni (GNAMPA) of the Istituto Nazionale di Alta Matematica (INdAM) and by the Universit\`a degli Studi di Catania, PIACERI 2020-2022, Linea di intervento 2, Progetto ”MAFANE”.
\bigskip
\bigskip
\bigskip
\bigskip
\centerline {\bf References}\par
\bigskip
\bigskip
\noindent
\noindent
[1]\hskip 5pt M. AIT MANSOUR, J. LAHRACHE and N. ZIANE, {\it Some applications of two minimax theorems}, Matematiche (Catania) {\bf 78} (2023), 405-414.\par
\smallskip
\noindent
[2]\hskip 5pt D. GIANDINOTO, {\it Further applications of two minimax theorems}, Matematiche (Catania) {\bf 77} (2022), 449-463.\par
\smallskip
\noindent
[3]\hskip 5pt B. RICCERI, {\it On the infimum of certain functionals}, in ``Essays in Mathematics and its Applications -
In Honor of Vladimir Arnold", Th. M. Rassias and P. M. Pardalos eds., 361-367, Springer, 2016.\par
\smallskip
\noindent
[4]\hskip 5pt B. RICCERI, {\it On a minimax theorem: an improvement, a new proof and an overview of its applications},
Minimax Theory Appl., {\bf 2} (2017), 99-152.\par
\smallskip
\noindent
[5]\hskip 5pt B. RICCERI, {\it Minimax theorems in a fully non-convex setting}, J. Nonlinear Var. Anal., {\bf 3} (2019), 45-52.\par
\smallskip
\noindent
[6]\hskip 5pt  M. SION {\it On general minimax theorems}, Pacific J. Math., {\bf 8} (1958), 171-176.\par
\smallskip
\noindent
[7]\hskip 5pt E. ZEIDLER, {\it Nonlinear functional analysis and its
applications}, vol. III, Springer-Verlag, 1985.\par

\bigskip
\bigskip
\bigskip
\bigskip
Department of Mathematics and Informatics\par
University of Catania\par
Viale A. Doria 6\par
95125 Catania, Italy\par
{\it e-mail address}: ricceri@dmi.unict.it

\bye